\documentclass[preprint,times,12pt]{article}

\usepackage{amsmath,amssymb,amsthm}
\usepackage{graphicx}
\usepackage{hyperref}
\usepackage{cleveref}
\usepackage{color}
\allowdisplaybreaks
\DeclareGraphicsExtensions{.pdf,.png,.jpg,.jpeg,.mps,.eps}

\usepackage{times}
\usepackage[all]{xypic}

\newtheorem{remark}{Remark}[section]
\newtheorem{theorem}{Theorem}[section]

\newtheorem{lemma}{Lemma}[section]                                                                                                              

\newtheorem{definition}{Definition}[section]

\numberwithin{equation}{section}

\def\R{\mathbb R}
\def\C{\mathbb C}

\def\Z{\mathbb Z}
\def\G{\mathsf{G}}
\def\H{\mathbb H}

\def\X{\mathcal X}

\def\P{{\mathcal P}}

\newcommand{\eps}{\varepsilon}
\def\tr{{\rm tr}}
\newcommand{\PSL}{{\rm PSL}}
\newcommand{\SL}{{\rm SL}}

\def\T{{\mathcal T}}

\begin{document}

\title{On 2-antiparallel encounters on factors of the hyperbolic plane}
\author{{\sc Huynh Minh Hien} \\[1ex] 
	Department of Mathematics,\\
	Quy Nhon University, Binh Dinh, Vietnam\\
	e-mail: huynhminhhien@qnu.edu.vn
	\date{} 
}
\maketitle

\begin{abstract} In this paper, we consider the geodesic flow on factors of the hyperbolic plane.  We prove that a periodic orbit including a 2-antiparallel encounter
	has a partner orbit. We construct the partner orbit and give an estimate for
	the action different between the orbit pair. Then we apply the result to reprove the accuracy of Sieber/Richter's prediction in [Sieber and Richter, 2001].
\end{abstract}

\section{Introduction}
The two-point correlator function of a classical dynamical system can be illustrated as a double sum over periodic orbits
\begin{equation}\label{formfactor}
K(\tau)=\Big\langle \frac{1}{T_H}\sum_{\gamma,\gamma'}A_\gamma 
A_{\gamma'}^* e^{\frac{i}{\hbar}(S_\gamma-S_{\gamma'})}
\delta\Big(\tau T_H-\frac{T_\gamma+T_{\gamma'}}{2}\Big) \Big\rangle,
\end{equation}
where $\langle\cdot\rangle$ abbreviates the  average over the energy and
over a small time window, $T_H$ denotes the Heisenberg time and $A_\gamma$, $S_\gamma$, and $T_\gamma$ 
are the amplitude, the action, and the period of the orbit $\gamma$, respectively. 
As one is interested in the semiclassical limit $\hbar \to 0$,  it is expected that only orbit pairs $\gamma, \gamma'$ such that $S_\gamma-S_{\gamma'}\sim \hbar$ or small. Formulated in more mathematical
terms, for a classical chaotic dynamical system the problem is to determine the periodic
orbit pairs $\gamma,\gamma'$ such that $S_\gamma$ is  close to $S_{\gamma'}$, and then calculate \eqref{formfactor}.

This was first considered by Sieber and Richter \cite{Sieber1,SieberRichter} 
who predicted that a given periodic orbit with a small-angle self-crossing in configuration space 
will admit a partner orbit with almost the same action. The original orbit and its partner are called  a Sieber-Richter pair. In phase space, a Sieber-Richter pair contains a region where two stretches of each  orbit are almost mutually time-reversed
and one addresses this region as a {\em $2$-encounter} or, more strictly, a {\em $2$-antiparallel encounter};
the `2' stands for two orbit stretches which are close in configuration space, 
and `antiparallel' means that the two stretches have opposite directions. 
The accuracy of Sieber/Richter's prediction was completely proven by Huynh/Kunze in \cite{HK}. In that paper the authors considered the geodesic flow on compact factors of the hyperbolic plane. It was shown in \cite{HK} that
a $T$-periodic orbit of the geodesic flow crossing itself in 
configuration space at a time $T_1$ has $9|\sin(\phi/2)|$-partner orbit and
the action difference between them is approximately equal
$\ln(1-(1+e^{-T_1})(1+e^{-(T-T_1)})\sin^2(\phi/2)))$ with the estimated error 
$12\sin^2(\phi/2)e^{-T}$, where $\phi$ is the crossing angle. 

Periodic orbits with $L$-parallel encounters was investigated by  M\"uller et al. in  \cite{HMBH,mueller2004,mueller2005}. 
 We speak of an {\em $L$-encounter} when $L$ stretches of a periodic orbit are mutually close to each other up to time reversal. 
 In other words, 
 all the $L$ stretches must intersect a small Poincar\'e section.
  M\"uller et al. used combinatorics to count the number of partner orbits and provided an approximation for the action difference, but a construction of partner orbits and an error bound of the approximation 
 had not been derived. Then, Huynh \cite{HMH2} continued considering the hyperbolic dynamical system in \cite{HK} to deal with the technically more involved higher-order encounters. The author proved that a given periodic
orbit including an $L$-parallel encounter has $(L - 1)! - 1$ partner orbits, constructed partner orbits and
gave estimates for the action differences between orbit pairs. Furthermore,   mathematical definitions for `encounters', `partner orbits', etc. were also arrived in \cite{HMH2}.

In the case of $L$-antiparallel encounter with general $L$, the problem is a very complicated and it is still open. In this paper, we only consider the problem for  $L=2$. We prove that a periodic orbit with $2$-antiparallel encounter has a partner orbit. If the space is compact, then the partner is unique. Then we apply the result to prove Sieber/Richter's prediction and derive a better estimate for the action difference. 
\smallskip

The paper is organized as follows. In Section 2 we introduce some necessary background materials. Then in Section 3 we consider periodic orbits with 2-antiparallel encounters. We prove that a periodic orbit with a 2-antiparallel encounter has a partner orbit.
Then we apply this result to reprove the accuracy of Sieber/Richter's prediction in \cite{SieberRichter}. 

\medskip

\noindent 

\setcounter{equation}{0}
\section{Preliminaries}

We consider the geodesic flow on factor $\Gamma\backslash \H^2$, 
where $\H^2=\{z=x+iy\in \C:\, y>0\}$ is the hyperbolic plane endowed  
with the hyperbolic metric $ds^2=\frac{dx^2+dy^2}{y^2}$ and
$\Gamma $ is a discrete subgroup of the projective Lie group $\PSL(2,\R)=\SL(2,\R)/\{\pm E_2\}$. The group $\PSL(2,\R)$  acts transitively on $\H^2$ by 
M\"obius transformations
$z\mapsto \frac{az+b}{cz+d}$.
If the action is free (of fixed points), then the factor $\Gamma\backslash\H^2$  has a Riemann surface structure.
Such a surface is a closed Riemann surface of genus at least $2$ 
and has the hyperbolic plane $\H^2$ as the universal covering. If the space $\Gamma\backslash\H^2$ is compact, then all the elements in $\Gamma\setminus\{e\}$
are hyperbolic, i.e., $\tr(\gamma)=|a+d|>2$ for $\gamma=\Big\{\pm \scriptsize\Big(\begin{array}{cc}a &b\\ c&d \end{array} \Big)\Big \}\in \Gamma\setminus\{e\}$. 
The geodesic flow $(\varphi_t^\X)_{t\in \R}$ on the unit tangent bundle $\X=T^1(\Gamma\backslash\H^2)$
goes along the unit speed geodesics on $\Gamma\backslash\H^2$. 
On the other hand, the unit tangent bundle $T^1(\Gamma\backslash\H ^2)$
is isometric to the quotient space 
$\Gamma\backslash \PSL(2,\R)=\{\Gamma g,g\in\PSL(2,\R)\}$, 
which is the system of right co-sets of $\Gamma$ in $\PSL(2,\R)$, by an isometry
$\Xi$.
Then  the geodesic flow $(\varphi_t^\X)_{t\in\R}$ can be equivalently expressed as the natural 
`quotient flow' $\varphi^X_t(\Gamma g)=\Gamma g a_t$ 
on $X=\Gamma\backslash\PSL(2,\R)$  associated to the flow $\varphi^\G_t(g)=g a_t$ on $\G:=\PSL(2,\R)$
by the conjugate relation 
\[\varphi_t^\X=\Xi^{-1}\circ\varphi^X_t\circ\Xi\quad \mbox{for all}\quad t\in\R.\]
Here $a_t\in\PSL(2,\R)$ denotes the equivalence class obtained from the matrix $A_t=\scriptsize\Big(\begin{array}{cc}
e^{t/2} & 0\\ 0 & e^{-t/2}
\end{array}\Big)\in\SL(2,\R)$. 

\smallskip 
There are some more advantages to work on $X=\Gamma\backslash\PSL(2,\R)$
rather than on $\X=T^1(\Gamma\backslash\H^2)$. One can calculate explicitly the stable and unstable manifolds 
at a point $x=\Gamma g\in X$ to be
\[W^s_X(x)=\{\Gamma gb_s,s\in\R\}
\quad \mbox{and}\quad W^u_X(x)=\{\Gamma gc_u, u\in\R\},\]
where  $b_s,c_u\in\PSL(2,\R)$ denote
the equivalence classes obtained from 
$B_s=\scriptsize\Big(\begin{array}{cc}1 &s\\ 0&1 \end{array} \Big), \ C_u=\scriptsize\Big(\begin{array}{cc}
1&0\\ u&1
\end{array}\Big)\in\SL(2,\R)$. If the space is compact, the flow $(\varphi^X_t)_{t\in\R}$
is hyperbolic.

General references for this section are \cite{bedkeanser,einsward,KatHas}, 
and these works may be consulted for the proofs to all results which are stated above.

For $\phi\in\R$, denote by $d_\phi\in\PSL(2,\R)$ the equivalence class obtained from $D_\phi=\scriptsize\Big(\begin{array}{cc}\cos(\phi/2) &-\sin(\phi/2)\\ \sin(\phi/2)&\cos(\phi/2) \end{array} \Big)\in \SL(2,\R)$. 
\begin{lemma}\label{adpi} (a) The following relations hold for $t\in\R$: 
\begin{equation}\label{dpi}a_t d_{\pi}=d_{\pi} a_{-t},\quad 
 b_t d_{\pi}=d_{\pi} c_{-t},\quad 
 c_t d_{\pi}=d_{\pi} b_{-t}.
 \end{equation} 
(b Let $g=[G]\in {\rm PSL}(2, \R)$ for $G=\scriptsize \Big(\begin{array}{cc} a & b \\ c & d\end{array}\Big)
\in {\rm SL}(2, \R)$. If $a\neq 0$, then $g=c_u b_s a_t$ for 
\begin{equation}\label{tsu-def}
 t=2\ln |a|,\quad s=ab,\quad u=\frac{c}{a}. 
\end{equation} 
\end{lemma}
\noindent {\bf Proof\,:} (a) In ${\rm SL}(2, \R)$ we calculate 
\begin{eqnarray*} 
	A_t D_{\pi} & = & \Bigg(\begin{array}{cc} e^{t/2} & 0 \\
		0 & e^{-t/2}\end{array}\Bigg)\Bigg(\begin{array}{cc} 0 & 1 \\
		-1 & 0\end{array}\Bigg)=\Bigg(\begin{array}{cc} 0 & e^{t/2} \\
		-e^{-t/2} & 0\end{array}\Bigg)=\Bigg(\begin{array}{cc} 0 & 1 \\
		-1 & 0\end{array}\Bigg)\Bigg(\begin{array}{cc} e^{-t/2} & 0 \\
		0 & e^{t/2}\end{array}\Bigg)
	\\ & = & D_{\pi} A_{-t} 
\end{eqnarray*} 
which upon projection yields the first one. The argument is analogous for the others.

\smallskip 
\noindent
(b)  Let $(t, s, u)$ be given by (\ref{tsu-def}). To begin with, 
\[C_u B_sA_t=\Bigg(\begin{array}{cc} e^{t/2} & 0 \\ 0 & e^{-t/2}\end{array}\Bigg)
\Bigg(\begin{array}{cc} 1 & 0 \\ u & 1\end{array}\Bigg)
\Bigg(\begin{array}{cc} 1 & s \\ 0 & 1\end{array}\Bigg)
=\Bigg(\begin{array}{cc} e^{t/2} & se^{-t/2} \\ ue^{t/2} & (1+su)e^{-t/2}\end{array}\Bigg). \] 
If $a>0$, then $e^{t/2}=a$, $se^{-t/2}=b$, $ue^{t/2}=c$, and $(1+su)e^{-t/2}=(1+bc)/a=d$, 
using that $ad-bc=1$. Thus $C_u B_sA_t=G$ and $ c_u b_sa_t=g$. 
If $a<0$, then $e^{t/2}=-a$, $se^{-t/2}=-b$, $ue^{t/2}=-c$, and $(1+su)e^{-t/2}=-(1+bc)/a=-d$, 
and hence $C_u B_sA_t=-G$ which yields once again that $c_u b_sa_t=g$. 
{\hfill$\Box$}\bigskip

\begin{figure}[ht]
	\begin{center}
		\begin{minipage}{0.8\linewidth}
			\centering
			\includegraphics[angle=0,width=1.1\linewidth]{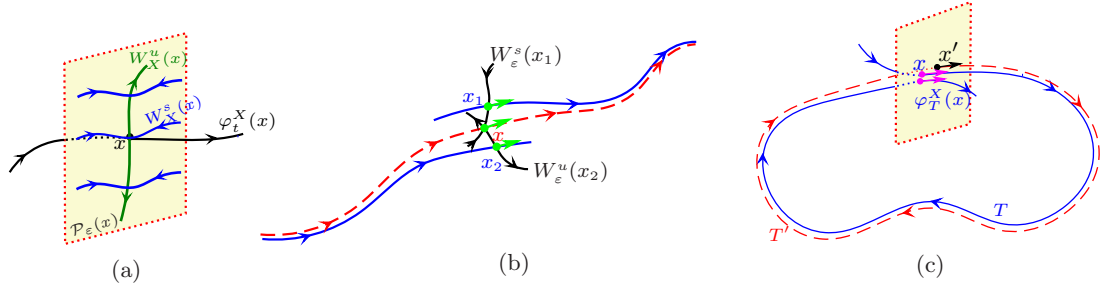}
		\end{minipage}
	\end{center}
	\caption{(a) Poincar\'e section \ \ (b) Shadowing lemma  \ \ (c)  Anosov closing lemma}\label{pss}
\end{figure}

We recall the definition of Poincar\'e sections, the shadowing lemma and the Anosov closing lemma from \cite{HK}.
\begin{definition}\label{Poindn1}\rm 
	Let $x\in X$ and $\eps>0$. The {\em Poincar\'e section} of radius $\eps$ at $x$ is 
	\begin{equation*}
	\P_\eps(x)=\{\Gamma(g c_ub_s) : |u|<\eps, |s| <\eps\},
	\end{equation*}
	where $g\in \G$ is such that $x=\Gamma g$ (see Figure \ref{pss}\,(a)).
\end{definition}
See Figure \ref{pss}\,(b)\&(c) for an illustration of the next two results. 
\begin{theorem}[Shadowing lemma]\label{shadlemII} 
	If $\eps>0$, $x_1, x_2\in X$, and $x\in W^s_{X,\,\eps}(x_1)\cap W^u_{X,\,\eps}(x_2)$, then 
	\[ d_X(\varphi^X_t(x_1), \varphi^X_t(x))<\eps e^{-t}
	\quad\mbox{for all}\quad t\in [0, \infty[ \]
	and 
	\[ d_X(\varphi^X_t(x_2), \varphi^X_t(x))<\eps e^t 
	\quad\mbox{for all}\quad t\in\, ]-\infty, 0]; \]
	recall $W^s_{X,\,\eps}(x)=\{\Gamma gb_s:\ |s|<\eps \}, W^u_{X,\,\eps}(x)=\{\Gamma gc_u: \ |u|<\eps \}$, for any $g\in\PSL(2,\R)$ such that $\Gamma g=x$.
\end{theorem}
\begin{theorem}[{Anosov closing lemma}]\label{anosov1}
	Suppose that $\eps\in\, ]0, \frac{1}{4}[$, 
	$x\in X$, $T\ge 1$, and $\varphi^X_T(x)\in {\cal P}_\eps(x)$. 
	Let $x=\Gamma g$ and  $\varphi^X_T(x)=\Gamma gc_ub_s$ 
	for $g\in\PSL(2,\R), |u|<\eps,|s|<\eps$. 
	Then there are $x'\in {\cal P}_{2\eps}(x)$
	and $T'\in\R$ so that 
	\begin{equation*}\label{x'}\varphi^X_{T'}(x')=x'\quad \mbox{and}\quad d_X(\varphi^X_t(x),\varphi^X_t(x'))<2 (|u|+|s|)\quad \mbox{for all}\quad t\in[0,T].
	\end{equation*}
	Furthermore, 
	\begin{equation}\label{TT'2}
	e^{T'/2}+e^{-T'/2}=e^{T/2}+e^{-T/2}+us e^{-T/2}
	\end{equation} 
and
	\begin{equation*}\label{T-T'1}
	\Big|\frac{T'-T}2-\ln(1+us)\Big| <5|us|e^{-T}.
	\end{equation*}
\end{theorem}
\begin{figure}[ht]
	\begin{center}
		\begin{minipage}{0.8\linewidth}
			\centering
			\includegraphics[angle=0,width=1.1\linewidth]{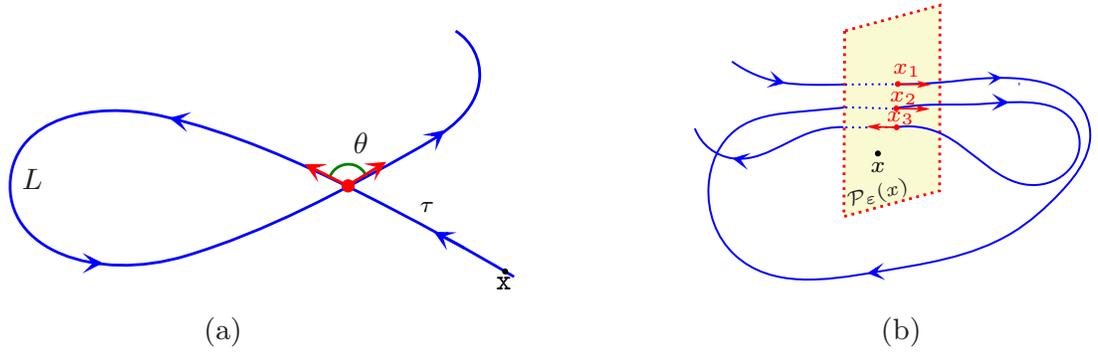}
		\end{minipage}
	\end{center}
	\caption{(a) An orbit with a self-crossing in configuration space \ \ \ (b) An orbit with a 3-antiparallel encounter}\label{se}
\end{figure}

See Figure \ref{se}\,(a) for an illustration for the next result.
\begin{theorem}[Self-crossings]\label{selfcthm} 
	Suppose that all elements of $\Gamma\setminus\{e\}$ are hyperbolic 
	and let  $\tau\in\R,\, L>0,\,\theta\in\,]0,\pi[$, and ${\mathtt x}\in \X$ be given. 
	The orbit of ${\mathtt x}$ under the geodesic flow 
	${(\varphi^\X_t)}_{t\in\R}$ crosses itself in configuration space at the 
	time $\tau$, at the angle $\theta$, and creates a loop of length $L$ if and only if
	\begin{equation}\label{sc1} \mbox{either}\quad \Gamma ga_{\tau+ L}=\Gamma g a_\tau d_\theta \quad\mbox{or}\quad
	\Gamma g_{\tau+L}=\Gamma g a_{\tau} d_{-\theta} 
	\end{equation}
	holds for any $g\in\PSL(2,\R), \Gamma g=\Xi({\mathtt x})$.
	Furthermore, 
	\begin{equation}\label{sc2}
	e^{-L}<\cos^2\Big(\frac{\theta}{2}\Big). 
	\end{equation}
\end{theorem}

The following definitions can be found in \cite{HMH2}. 
\begin{definition}[Time reversal]\label{tr}\rm
	The {\em time reversal map} $\T:X \rightarrow X$ is defined by 
	\begin{equation*}\label{revedn}
	\T(x)=\Gamma gd_\pi \quad \mbox{for}\quad  x=\Gamma g\in X,
	\end{equation*}
	where $d_\pi\in\PSL(2,\R)$ is the equivalence class of the matrix 
	$D_\pi=\scriptsize\Big(\begin{array}{cc} 0&1\\-1&0\end{array}\Big)\in{\rm SL(2,\R)}$.
\end{definition} 
Using Lemma \ref{adpi}\,(a), we have
\begin{equation}\label{tr2}
\varphi^X_t(\T(x))=\T(\varphi^X_{-t}(x))
\quad \mbox{for}\quad x\in X\quad\mbox{and}\quad t\in\R.
\end{equation}

\begin{definition}[Orbit pair/Partner orbit]\label{partnerdf1}\rm 
	Let $\eps>0$ be given. 
	Two
	given $T$-periodic orbit $c$ and $T'$-periodic orbit $c'$ of the flow
	$(\varphi^X_t)_{t\in\R}$ are called an {\em $\eps$-orbit pair}
	if 
	there are $L\geq 2, L\in \Z$ and two decompositions  of $[0,T]$ and $[0,T']:$ $0=t_0<\cdots<t_L=T$
	and $0=t_0'<\cdots<t_L'=T'$, 
	and a permutation $\sigma:\{0,1,\dots,L-1\}\rightarrow 
	\{0,1,\dots, L-1\}$ such that 
	for each $j\in\{0,\dots,L-1\}$, 
	either
	\begin{equation*}
	d_X(\varphi^X_{t+t_j}(x), \varphi^X_{t+t_{\sigma(j)}'}(x'))<\eps
	\quad \mbox{for all}\quad t\in [0,t_{j+1}-t_j]
	\end{equation*}
	or
	\begin{equation*}
	d_X\Big(\varphi^X_{t+t_j}(x), \varphi^X_{t-t_{\sigma(j)+1}'}(\T(x'))\big)<\eps
	\quad \mbox{for all}\quad t\in [0,t_{j+1}-t_j]
	\end{equation*}
	holds for some $x\in c$ and $x'\in c'$.
	Then $c'$ is called an {\em $\eps$-partner orbit} of $c$ and vice versa. 
\end{definition}
Roughly speaking, two periodic orbits are an $\eps$-orbit pair if they are $\eps$-close to each other in configuration space, not for the whole time, since otherwise they would be identical, but they decompose to the same number of parts and any part of one orbit is $\eps$-close to some part of the other.

\begin{definition}[Encounter]\rm Let $\eps>0$ and $L\in \Z, L\geq 2$ be given.
	We say that a periodic orbit $c$ has an {\em$(L,\eps)$-encounter} if
	there are $x\in X$, $x_1,\dots, x_L\in c$ such that for each $j\in\{1,\dots,L\}$,
	\[\mbox{either}\quad x_j\in \P_\eps(x)\quad \mbox{or}\quad \T(x_j)\in \P_\eps(x).\]
	If either $x_j\in\P_\eps(x)$ holds for all $i=1,\dots,L$ or $\T(x_j)\in\P_\eps(x)$ holds for all $j=1,\dots,L$  then
	the encounter is called {\em parallel encounter}; otherwise it is called {\em antiparallel encounter} (see Figure \ref{se}\,(b)).
\end{definition}

\section{Main results}

\subsection{2-Antiparallel encounters} 
In this section we only consider 2-antiparallel encounters.
It is impossible to reconnect the ports in 2-parallel encounter to get a new (genuine) partner orbit; but in the case of antiparallel encounter, we have the following result.
\begin{theorem}\label{2anti} Suppose that all the elements in $\Gamma\setminus\{e\}$ are hyperbolic and let $\eps>0$. 
If a periodic orbit $c$ of the flow $(\varphi^X_t)_{t\in\R}$ on $X$
	with period $T>1$ has a $(2,\eps)$-antiparallel encounter, then it has a partner. Furthermore,
	let $x,y\in c$, $x=\Gamma g$ and $\T(y)=\Gamma g c_ub_s\in\P_\eps(x)$ for $g\in\G, |u|<\eps,|s|<\eps$; $\varphi^X_{T_1}(x)=y, 0< T_1<T$. Then the partner is $\eps'$-partner with $\eps'=\eps+2(|u-se^{-T_1}|+|s-ue^{T_1-T}|)<9\eps$ and the action
	difference between the orbit pair satisfies
	\begin{equation}\label{TT'}
	\Big|\frac{T'-T}{2}-\ln\big(1+(u-se^{-T_1})(s-ue^{T_1-T})\big)\Big|\leq |(u-se^{-T_1})(s-ue^{T_1-T})|e^{-T},
	\end{equation}
	where $T'$ is the period of the partner. 
\end{theorem}
\noindent{\bf Proof\,:}
Let $x=\Gamma g$ and write $y=\Gamma h$ with  $g,h\in \G$ and set $g'=gd_\pi, h'=hd_\pi, T_2=T-T_1$. Then by the assumption
${\cal T}(y)=\Gamma h'=\Gamma g c_ub_s$ or $\Gamma h=\Gamma g'b_{-u}c_{-s}$ due to Lemma \ref{adpi}\,(a).
This implies that $w:=\Gamma h' b_{-s}=\Gamma g c_u\in W_{X,\,\eps}^s(y')\cap W_{X,\,\eps}^u(x)$.
By the shadowing lemma (Theorem \ref{shadlemII}),
\begin{equation}\label{t>0'}
d_X(\varphi^X_t(y'),\varphi^X_t(w))<\eps e^{-t}\quad \mbox{for all}\quad t\in[0,\infty[
\end{equation}
and
\begin{equation}\label{t<0'}
d_X(\varphi^X_t(x),\varphi^X_t(w))<\eps e^{t}\quad \mbox{for all}\quad t\in\,]-\infty,0].
\end{equation}
Putting $\hat w=\varphi^X_{-T_2}(w)=\Gamma g c_u a_{-T_2}=\Gamma h'b_{-s}a_{-T_2}$, we claim that $\varphi^X_T(\hat w)\in \P_{2\eps}(\hat w)$. Indeed,
\begin{eqnarray*}
	\varphi^X_T(\hat w)&=&\Gamma h' b_{-s}a_{T_1}
	=\Gamma g'a_{-T_1}b_{-s}a_{T_1}
	=\Gamma g' b_{-se^{-T_1}}
	=\Gamma h c_sb_ub_{-se^{-T_1}}
	\\
	&=& \Gamma ga_{-T_2}c_sb_{u-se^{-T_1}}
	=\Gamma(gc_ua_{-T_2})(a_{T_2}c_{-u}a_{-T_2}c_sb_{u-se^{-T_1}})
	\\
	&=& \Gamma gc_ua_{-T_2}(c_{-ue^{-T_2}}c_sb_{u-se^{-T_1}})
	=\Gamma (g c_ua_{-T_2})(c_{-ue^{-T_2}+s}b_{u-se^{-T_1}})
	\\
	&=& \Gamma (gc_ua_{-T_2}) (c_{u'} b_{s'})
\end{eqnarray*}
with
\begin{equation}\label{u's'}
u'=s-ue^{-T_2},\quad s'=u-s e^{-T_1}.
\end{equation}
Apply the Anosov closing lemma (Theorem \ref{anosov1}) to obtain $v\in X, T'\in\R$ such that
$\varphi^X_{T'}(v)=v$,\[\quad \Big|\frac{T'-T}{2}-\ln(1+u's')\Big|<5|u's'|e^{-T},\] and
\begin{equation}\label{vw'}
d_X(\varphi^X_t(\hat w),\varphi^X_t(v))<2(|u'|+|s'|)\quad\mbox{for all}\quad t\in [0,T].
\end{equation}
For $t\in[0,T_1]$, it follows from \eqref{t<0'} and \eqref{vw'} that
\begin{eqnarray*}
	d_X(\varphi^X_t(v),\varphi^X_t(y'))
	&\leq& d_X(\varphi^X_t(v),\varphi^X_t(\hat w))
	+d_X(\varphi^X_t(\hat w),\varphi^X_t(y'))\\
	&\leq&d_X(\varphi^X_t(v),\varphi^X_t(\hat w))
	+d_X(\varphi^X_t(\varphi^X_{-T_1}(w)),\varphi^X_t(\varphi^X_{-T_1}(x')))\\
	&\leq& 2(|u'|+|s'|)+d_X(\varphi^X_{t-T_1}(w),\varphi^X_{t-T_1}(x'))\\
	&<&2(|u'|+|s'|)+\eps=\eps'.
\end{eqnarray*}
Similarly, for $t\in[T_1,T]$, it follows from \eqref{t>0'} and \eqref{vw'} that
\begin{eqnarray*}
	d_X(\varphi^X_t(v),\varphi^X_t(y))
	&=& d_X(\varphi^X_t(v),\varphi^X_{t-T_1}(y))
	\leq d_X(\varphi^X_t(v),\varphi^X_{t}(\hat w))+d_X(\varphi^X_{t}(\hat w),\varphi^X_{t-T_2}(y))\\
	&\leq& 2(|u'|+|s'|)+ d_X(\varphi^X_{t-T_2}(w),\varphi^X_{t-T_2}(y))\\
	&\leq& 2(|u'|+|s'|)+\eps=\eps'.
\end{eqnarray*}
We can easily check that the partner is a $\eps'$-partner orbit in the sense of Definition \ref{partnerdf1}. 
{\hfill$\Box$}
\begin{remark}\label{TT'rm}
	It follows from \eqref{TT'2} and \eqref{u's'} that 
	\begin{itemize}
		\item[(i)] $T'>T$ if and only if $ (s-ue^{-T_2})(u-s e^{-T_1})>0$;
		\item[(ii)] $T'<T$ if and only if $ (s-ue^{-T_2})(u-s e^{-T_1})>0$;
		\item[(iii)] $T'=T$ if and only if $(s-ue^{-T_2})(u-s e^{-T_1})=0$.
	\end{itemize}
\end{remark}
{\hfill$\Diamond$}
\subsection{An application to Sieber-Richter pairs}

\begin{figure}[ht]
	\begin{center}
		\begin{minipage}{0.6\linewidth}
			\centering
			\includegraphics[angle=0,width=1.1\linewidth]{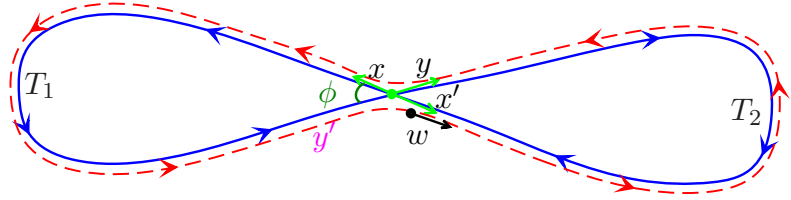}
		\end{minipage}
	\end{center}
	\caption{Small self-crossing angle and partner orbit}\label{mainthm}
\end{figure}

A periodic orbit with a small-angle self-crossing has 2 almost mutually time-revered stretches. This means that the orbit crosses the Poincar\'e section of a point in this orbit and Theorem \ref{2anti} may be applied. 
\begin{theorem}\label{existthm}
	If a periodic orbit of the geodesic flow ${(\varphi^\X_t)}_{t\in\R}$ on $\X=T^1(\Gamma\backslash\H^2)$
	with the period $T\geq 1$ crosses itself in configuration space at a time $T_1\in\,]0,T[$
	and at an angle $\theta$ such that $0<\phi<\frac13$
	for $\phi=\pi-\theta$, 
	then it has  a $6|\sin(\phi/2)|$-partner orbit. 
	Furthermore, $T'<T$ for the period of the partner and the action difference
	satisfies
	\begin{eqnarray}\label{actdiff}
	\Big|\frac{T'-T}{2}-\ln\big(1-\sin^2(\phi/2)({\cos^{-2}(\phi/2)}+e^{-T_1})(\cos^2(\phi/2)+e^{T_1-T})\big)\Big|\leq 2\sin^2(\phi/2)e^{-T}.
	\end{eqnarray}
\end{theorem}
\noindent {\bf Proof\,:}  Let the orbit of ${\mathtt x}\in {\cal X}=T^1(\Gamma\backslash\H^2)$ 
be $T$-periodic ($T$ is the prime period)
and such that it has a self-crossing of angle $\theta$ 
in configuration space at the time $T_1\in \,]0, T[$, i.e., we have 
\begin{equation}\label{seya}
\varphi_{T_1}^{{\cal X}}({\mathtt x})={\mathtt y},
\quad\varphi_{T_2}^{{\cal X}}({\mathtt y})={\mathtt x},
\end{equation} 
where $T=T_1+T_2$; see Figure \ref{mainthm}.

In addition, assume that $|\phi|<\frac13$ with $\phi=\pi-\theta$. Then in particular 
\begin{equation}\label{sinphi/2} 
\Big|\sin\Big(\frac{\phi}{2}\Big)\Big|\le\frac{|\phi|}{2}
<\frac16
\end{equation} 
holds.
Set $x=\Xi({\mathtt x})$
and $y=\Xi ({\mathtt y})$. This follows from Theorem  \ref{selfcthm} that
\begin{equation}
\Gamma g=\Gamma h d_\theta \quad\mbox{or}\quad \Gamma g=\Gamma hd_{-\theta}
\end{equation}
with some $h,g\in\PSL(2,\R)$ such that $\Gamma g=x$ and $\Gamma h=y$.  We only consider the first case, the latter is similar. Write 
$x'={\cal T}(x)$ and $y'={\cal T} (y)$; recall the notation ${\mathcal T}$ from Definition \ref{tr}. Then
\[\Gamma gd_{\pi-\theta}=\Gamma h d_\pi=y'\quad \mbox{or}\quad \Gamma g d_\phi=y'.\]
We write 
\begin{equation}\label{dexpress} 
d_{\phi}=c_ub_s a_\tau,
\end{equation} 
where 
\[ \tau=2\ln(\cos(\phi/2)),\quad u=\tan(\phi/2),\quad s=-\sin(\phi/2)\cos(\phi/2). \]  
By (\ref{sinphi/2}), we have 
\begin{equation}\label{cosphi/2}
\cos\Big(\frac{\phi}{2}\Big)>\frac{5}{6}.
\end{equation} 
Then 
\begin{equation}\label{tausu-bd}  
|u|=|\tan(\phi/2)|\le \frac65|\sin(\phi/2)|=:\eps,
\ |s|=|\sin(\phi/2)\cos(\phi/2)|\le |\sin(\phi/2)|<\eps,
\end{equation} 
and 
\[ |\tau|=|\ln(1-\sin^2(\phi/2))|\le 2\sin^2(\phi/2)
\le \frac12\, \eps^2, \]
due to $|\ln(1+z)|\le 2|z|$ for $|z|\le 1/2$. Denote 
$\tilde{y}=\varphi^X_{-\tau}(y)$. This leads to
\[{\cal T}(\tilde y)=\varphi^X_\tau(y')=\Gamma h a_{-\tau}=\Gamma g c_ub_s\in \P_\eps(x),\]
using \eqref{tr2}. We are in a position to apply Theorem \ref{2anti}  to have $v\in X$ and
$T'\in \R$ such that $\varphi^X_{T'}(v)=v$ and the orbit of $v$ is $\eps'$-partner of the orbit of $x$, where $\eps'$ is determined later. 
Note that 
\begin{eqnarray*}
	(s-ue^{-T_2})(u-s e^{-T_1})&=&\big(-\sin(\phi/2)\cos(\phi/2)-\tan(\phi/2)e^{-T_2}\big)\big(\tan(\phi/2)+\sin(\phi/2)
	\cos(\phi/2)e^{-T_1}\big)\\
	&=& -\sin^2(\phi/2)(\cos^{-2}(\phi/2)+e^{-T_1})(\cos^2(\phi/2)+e^{-T_2})<0
\end{eqnarray*}
implies $T'<T$ owing to Remark \ref{TT'rm}. 
Furthermore, using \eqref{sc2} and \eqref{sinphi/2}
\begin{equation}\label{T122} e^{-T_1}<\sin^2(\phi/2)<\frac{1}{36}\quad\mbox{ and}\quad e^{-T_2}<\sin^2(\phi/2)<\frac{1}{36}, 
\end{equation}
yields \[|(s-ue^{-T_2})(u-s e^{-T_1})|<2\sin^2(\phi/2)e^{-T}\]
and hence
\begin{eqnarray*}
\Big|\frac{T-T'}{2}-\ln\big(1-\sin^2(\phi/2)({\cos^{-2}(\phi/2)}+e^{-T_1})(\cos^2(\phi/2)+e^{-T_2})\big)\Big|
 \leq 2 \sin^2(\phi/2)e^{-T}
\end{eqnarray*} 
which is \ref{actdiff}.
Finally, 
\[\eps'=\eps+2(|s-ue^{-T_2}|+|u-s e^{-T_1}|)\leq \frac{6}{5}|\sin(\phi/2)|+\frac{68}{15}|\sin(\phi/2|)<6|\sin(\phi/2)|\]
by \eqref{tausu-bd} and \eqref{T122}; and therefore the partner is a $6|\sin(\phi/2)|$-partner orbit. 
{\hfill$\Box$}

\begin{remark}\rm
	(i) Recall from \cite{HK} that the partner is a $9|\sin(\phi/2)|$-partner and the action difference satisfies
	\begin{eqnarray*}\label{act}
	\Big|\frac{T'-T}{2}-\ln\big(1-(1+e^{-T_1})(1+e^{-(T-T_1)})\sin^2(\phi/2)\big)\Big|\leq 12\sin^2(\phi/2)e^{-T}.
	\end{eqnarray*}
	This means that the orbits in a Sieber-Richter pair are estimated closer in this paper and   the estimate of the action difference is also better. 
	
	\noindent 
	 (ii) As mentioned in \cite{muellerthesis}, the partner orbit is avoided crossing. This means that it does not cross itself in encounter area. Conversely, a periodic orbit with 2 stretches almost mutually time-reversed and avoiding crossing has a partner orbit with a small-angle self-crossing in encounter area. Indeed, since the two stretches are almost mutually time-reversed, there are $x=\Gamma g$ and $y$ on that orbit such that $x$ and ${\cal T}(y)$ are very close. Using Lemma \ref{adpi} we write ${\cal T}(y)=\Gamma gc_ub_sa_\tau$ with $|u|,|s|<\eps$ for some small $\eps$.
Then $\varphi_{-\tau}^X({\cal T}(y))=\Gamma gc_ub_s$
implies that ${\cal T}(z)\in \P_\eps(x)$ for $z=\varphi_\tau^X(y)$ and we can apply 
Theorem \ref{2anti} to obtain a partner orbit for the original one. 

\smallskip 
\noindent
(iii) According to \cite{HK}, if the space $X=\Gamma\setminus\PSL(2,\R)$ is compact (or equivalently, $\Gamma\backslash\H^2$ is compact) and the crossing angle is small enough then the partner is unique.

\end{remark}

\end{document}